\documentclass[a4paper,10pt]{amsart}
\usepackage{amssymb,latexsym,graphicx,amsxtra,amsmath,amsthm}
\usepackage{hyperref}

\newcommand{\al}{\alpha}

\newcommand{\q}{g_n(\zeta)=\rho_n^{-k}[f_n(z_n+\rho_n\zeta)-a] }

\newcommand{\fr}{\mathcal F}
\newcommand{\ity}{\infty}
\newcommand{\C}{\mathbb{C}}

\numberwithin{equation}{section}
\newtheorem{theorem}{Theorem}[section]
\newtheorem{lemma}[theorem]{Lemma}

\theoremstyle{remark}

\thanks {The research work of the first author is supported by junior research fellowship from UGC India. }

\begin{document}

\title[NORMALITY AND SHARING VALUES ]{NORMALITY AND SHARING VALUES}

\author[G.Datt]{Gopal Datt}
\address{Department of Mathematics, University of Delhi,
Delhi--110 007, India}

\email{ggopal.datt@gmail.com,datt.gopal@ymail.com }

\author[S. Kumar]{Sanjay Kumar}

\address{Department of Mathematics, Deen Dayal Upadhyaya College, University of Delhi,
Delhi--110 007, India }

\email{sanjpant@gmail.com}

\begin{abstract}{In this paper, we obtain some normality criteria for families of holomorphic functions. these generalize some results of Fang, Xu, Chen and Hua.
}
\end{abstract}



\maketitle

\section{Introduction}\label{sec1}
We denote the  complex plane by $\C$, and the unit disk by $\Delta$. Let $f$ be a meromorphic function in$\C$. We say that $f$ is a normal function if there exits a positive $M$ such that$ f^{\#} (z) \leq M$ for all $z\in \C$, where $f^{\#}=\frac{|f'(z)|}{1+|f(z)|^2}$ denotes the spherical derivative of $f.$ \\

   A family $\fr$ of analytic functions on a domain $\Omega\subseteq\C$ is normal in $\Omega$ if every sequence of functions ${f_n}\subseteq\fr$ contains either a subsequence which converges to a limit function $f\not\equiv\ity$ uniformly on each compact subset of $\Omega$, or a subsequence which converges uniformly to $\ity$ on each compact subset.\\

In this paper, we use the following standard notation of value distribution theory,\\

\qquad\qquad\qquad\qquad\qquad $T(r,f); m(r,f); N(r,f); \overline{N(r,f)},\ldots$\\
We denote $S(r,f)$ any function satisfying\\

\qquad\qquad\qquad\qquad\qquad $S(r,f)=o\{T(r,f)\}$,  as $r\rightarrow +\ity,$\\
 possibly outside of a set with finite measure.\\

According to Bloch's principle every condition which reduces a meromorphic function in the plane to a constant, makes the family of meromorphic functions in a domain G normal. Rubel gave four counter examples to Bloch principle. \\Let $f$ and $g$ be two meromorphic functions in a domain $D \text{ and}\ a\in \C$. If $f-a$ and $g-a$ have the same number of zeros in $D$ (ignoring multiplicity). Then we say that $f$ and $g$ share the value $z=a$ IM.\\ Let us  recall the following known results that establish connection between shared values and normality. \\ Mues and Steinmetz proved the following result.\\

\begin{theorem}~\cite{MS}
{ Let $f$ be a non constant meromorphic  function in the plane. If $f$ and $f'$ share three distinct complex numbers $a_1$, $a_2$, $a_3$ then $f$ $\equiv$ $f'$.}

\end{theorem}
Wilhelm Schwick seems to have been the  first to draw a connection between normality and shared values. He proved the following theorem\\

\begin{theorem}
~\cite{Sch}{Let $\mathcal F$ be a family of meromorphic functions on a domain $G$ and $a_1$, $a_2$, $a_3$ be distinct complex numbers . If $f$ and $f'$ share $a_1$, $a_2$, $a_3$ for every $f\in \mathcal F$, then  $\mathcal F$ is normal in $G$.}
\end{theorem}

Chen and Hua proved the following theorem\\

\begin{theorem}~\cite{CH}
Let $\mathcal F$ be a family of holomorphic functions in a domain $D$. Suppose that there exists a non zero  $a\in \C$ such that for each function $f \in F$; $f, f'$ and $f''$ share the value $z=a$ IM in $D$. Then the family $\mathcal F$ is normal in $D$.
\end{theorem}

Fang and Xu improved their results by proving the following theorems\\

\begin{theorem}~\cite{FX}\label{th1.4}
Let $\mathcal F$ be a family of holomorphic functions on a domain $D$ and let $a$, $b$ be two distinct finite complex numbers such that $b \neq$ 0. If for any  $f \in \mathcal F$, $f$ and $f'$ share $z=a$ IM and $f(z)=b$ whenever $f'(z)=b$ then $\mathcal F$ is normal in $D$.\\
\end{theorem}

 \begin{theorem}~\cite{FX}\label{th1.5}
 {Let$\fr$ be a family of holomorphic functions in a domain $D$, and let $a$ be a non zero finite complex number. If for any $f\in \fr$ $f$ and $f'$ share $z=a$ IM and $f^{(k)}(z)=a, f^{(k+1)}(z)=a$ whenever $f(z)=a$. Then $\fr$ is normal in $D$.}
\end{theorem}

Finally,  Fang proved the following. \\

\begin{theorem}~\cite{MLF}\label{th1.6}
{Let $\fr$ be a family of meromorphic functions in a domain $D$ and let $a(z)$ be a non vanishing analytic function in $D$. If, for every function $f\in\fr, \ f$ and  $f'$  have the same zeros, and $f(z)=a(z)$ whenever $f'(z)=a(z),$ then $\fr$ is normal in $D.$}
\end{theorem}
More recently Xia and Xu improved theorem \ref{th1.6} by showing the following :\\

\begin{theorem}~\cite{XX}
{Let $\fr$ be a family of meromorphic functions in a domain  D, and k be a positive  integer, and let  $\varphi(z)(\not\equiv 0, \ity)$ be a non vanishing meromorphic function in D such that f and $ \varphi(z)$ have no common zeros for all $ f\in\fr$ and  $ \varphi(z)$   has no simple zeros in D, and all poles of $\varphi(z)$  have multiplicity at most k. If, for each $ f\in\fr,$ \\ $(1)$  all zeros of f have multiplicity at least $k+1$\\ $(2)$ $  f(z)=0 \ \text{whenever}\ f^{(k)}(z)=0$ and $ f(z)=\varphi(z)$ whenever $f^{(k)}(z)=\varphi(z)$\\ then $\fr$  is normal in D.}
\end{theorem}
\section{Main Theorems and  Lemmas}
We improve Theorem \ref{th1.4} and Theorem \ref{th1.5} by showing
the following.
\begin{theorem} Let $\mathcal F$ be a family of holomorphic functions on a  domain  $D$ such that all zeros of $f\in \fr$ are of multiplicity at least $k$, where $k$ is a positive integer.   Let $a$, $b$ be two distinct finite complex numbers such that $b \neq$ 0.  Suppose  for any  $f \in \mathcal F$ satisfies the following\\$(1)$: \  $f$ and $f^{(k)}$ share $z=a$ IM \\$(2)$: \ $f(z)=b$ whenever $f^{(k)}(z)=b$\\ then $\mathcal F$ is normal in $D$.
\end{theorem}
One may ask whether we can replace the values $a$ and $b$ by holomorphic functions.  We show in the following theorem that this is indeed the case.\\

\begin{theorem}
Let $\mathcal F$ be a family of holomorphic functions on a  domain  $D$ such that all zeros of $f\in \fr$ are of multiplicity at least $k$, where $k$ is a positive integer. Let $a(z),\ b(z),\ \al_o(z), \ \al_1(z) $ be holomorphic functions in $D$,  with $\al_0(z)\neq 0.$  If, for each $f\in \fr$, \\$(1)$: \ $ \  b(z)\neq 0$ \\ $(2)$: \    $ a(z)\neq b(z)$, and  {$b(z)-\al_1(z)a(z)-\al_0(z)a^{(k)}(z)\neq 0$}.\\
$(3)$: \ {$f(z)=a(z)\ \text{ if and only if}\ \ \al_0(z) f^{(k)}(z)+\al_1(z) f(z)=a(z)$}\\
$(4)$: \ {$ f(z)=b(z)\ \text {whenever}\ \ \al_0(z) f^{(k)}(z)+\al_1(z) f(z)=b(z)$ }\\
then $\fr$ is normal in $D.$

\end{theorem}
\textbf{Remark \ 1:} The hypothesis $a(z)\neq b(z)$ and {$b(z)-\al_1(z)a(z)-\al_0(z)a^{(k)}(z)\neq 0$} can not be dropped in Theorem 2.2.\\

\textbf{Example \ 1:} Let $D = \Delta = \{z: \mid z \mid < 1\}$ and $a(z) = b(z) = z^{k-1}$, $\al_o(z)= 1, \ \al_1(z)=0$ and \\

\qquad\qquad\qquad\qquad\qquad $\fr = \{e^{nz} - \frac{z^{k-1}}{n^k} + z^{k-1} : n = 1, 2, \ldots\}$.\\

Then for any $f\in \fr$, and \\

 \qquad\qquad\qquad\qquad\qquad $f = e^{nz} -\frac{z^{k-1}}{n^k} + z^{k-1} , \quad f^{(k)} = n^{k}e^{nz}$\\

 Clearly, conditions of Theorem 2.2 are satisfied. However, $\fr$ is not normal in $\Delta$.\\

This example  confirms that  $b(z)\neq 0$ is necessary in Theorem 2.2 as $f^{(k)}(z)\neq 0$.\\

\textbf{Example \ 2:} Let $D = \Delta =  \{z: \mid z \mid < 1\}$,  $k$ be a positive integer, $b(z)= b \ \text{(a non zero constant)} \ \text{ and}\  a(z) =((-1)^{k+1}+1)b$   and
\begin{equation}\notag
\fr=\{b\frac{(z-\frac{1}{n})^k}{k!}+\frac{(-1)^{k+1}}{k!(z-\frac{1}{n})}+a : n=1,2,\ldots\}
\end{equation}

Then, for every $f_n(z)\in \fr$,

\begin{equation}\notag
f_n(z)=b\frac{(z-\frac{1}{n})^k}{k!}+\frac{(-1)^{k+1}}{k!(z-\frac{1}{n})}+a,\qquad f_n^{(k)}(z)=b-\frac{1}{(z-\frac{1}{n})^{k+1}}
\end{equation}

Clearly, $f_n$ and $f_n^{(k)}$ share $a$ and $f_n^{(k)}(z)\neq b$,
so that $f_n(z)= b$ whenever  $f_n^{(k)}(z)=b$. But $\fr$ is not normal in $D$.

 \begin{theorem}Let $\fr$ be a family of holomorphic functions in a domain $D$ such that all zeros of $f\in \fr$ are of multiplicity at least $k$, where $k$ is a positive integer and let a be a non zero finite complex number. If for any $f\in \fr$ $f$ and $f^{(k)}$ share $z=a$ IM and $f^{(k+1)}(z)=a$ whenever $f(z)=a$. Then $\fr$ is normal in $D$.
\end{theorem}
We will use the tools of Fang and Xu which they used in their paper. For this we need  the following.\\

\begin{lemma}~\cite{Zalc}~\cite{Zalc 1}\label{lem1}(Zalcman's lemma)\\Let $\mathcal F$ be a family of holomorphic functions in the unit disk  $\Delta$ with the property that for every function $f\in \mathcal F$ , the zeros of $f$ are of multiplicity at least k. If $\mathcal F$ is not normal at $z_0$ in $\Delta$, then for 0 $\leq \alpha <k$, there exist\\
(a) a sequence of complex numbers $z_n \rightarrow z_0$, $|z_n|<r<1$\\
(b) a sequence of functions $f_n\in \mathcal F$ and\\
(c) a sequence of positive numbers $\rho_n \rightarrow 0$\\
such that $g_n(\zeta)=\rho_n^{-\alpha}f_n(z_n+\rho_n\zeta) $ converges to a non constant entire function $g$ on $\C$. Moreover $g$ is of order at most one . If $\mathcal F$ possesses the additional property that there exists $M > 0$ such that $|f^{(k)}(z)|\leq M$ whenever $f(z)=0$ for any $f \in \mathcal F$, then we can take $\alpha=k$.\\

\end{lemma}

\begin{lemma}~\cite{Hay}~\cite{Yang}\label{lem2}
{Let f be a non constant meromorphic function. Then for $ k \geq 1, \  b \neq 0,\infty,$\\

\quad \quad \quad  $T(r,f)\leq \overline N(r,f)+N(r,\frac{1}{f})+N(r,\frac{1}{f^{(k)}-b})-N(r,\frac{1}{f^{(k+1)}})+S(r,f)$}\\

\end{lemma}

\section{Proof of Theorem 2.1}

\begin{proof}Since normality is a local property, we assume that $D=\Delta=\{z:|z|<1\}$
Suppose,  $\fr$ is not normal in $D$; without loss of generality we assume that $\fr$ is not normal at the point $z_0$ in $\Delta$.  Then by Lemma \ref{lem1}, there exist \\(a) a sequence of complex numbers $z_n \rightarrow z_0$, $|z_n|<r<1$\\
(b) a sequence of functions $f_n\in \mathcal F$ and\\
(c) a sequence of positive numbers $\rho_n \rightarrow 0$\\
such that $\q$\\ converges locally uniformly to a non constant entire function $g$. Moreover $g$ is of order at most one.\\

Now we claim that   $g=0\text{ if and only if}\  g^{(k)}=a$ and $g^{(k)}\neq b$\\ Suppose,  $g(\zeta_0)=0$. then by Hurwitz's theorem there exist $\zeta_n;   \zeta_n\rightarrow \zeta_0$ such that\\

\qquad \qquad \qquad \qquad $g_n(\zeta_n)=\rho_n^{-k}[f_n(z_n+\rho_n\zeta_n)-a]=0 $\\

 Thus $f_n(z_n+\rho_n\zeta_n)=a$. Since $f_n$ and $ f_n^{(k)}$ share $z=a$ IM , we have\\

\qquad \qquad \qquad \qquad \qquad \qquad  $g_n^{k}(\zeta_n)=f_n^{(k)}(z_n+\rho_n\zeta_n)=a$\\ Hence $$g^{(k)}(\zeta_0)=\lim_{n \rightarrow \infty}g_n^{(k)}(\zeta_n)=a$$\\
Thus we have proved that $g^{(k)}=a \text{ whenever} \ g=0.$\\

\qquad   On the other hand, if $g^{(k)}(\zeta_0)=a$,  then there exist $\zeta_n; \zeta_n \rightarrow \zeta_0$ such that\\ $g_n^{(k)}(\zeta_n)=f_n^{(k)}(z_n+\rho_n\zeta_n)=a;\  n=1,2,\ldots$ hence $f_n(z_n+\rho_n\zeta_n)=a$ and $g_n(\zeta_n)=0$ for n=1,2,\ldots  thus $$g(\zeta_0)=\lim_{n\rightarrow \infty}g_n(\zeta_n)=0$$ this shows that $g=0$ whenever
$g^{(k)}=a$\\ Hence $g=0 \ \text{ if and only if }\  g^{k}=a$.\\

\qquad Next, we prove $g^{(k)}(\zeta)\neq b.$
 Suppose that there exist $\zeta_0$ satisfying $g^{(k)}(\zeta_0)= b$. Then, by Hurwitz's theorem, there exist a sequence $\zeta_n \rightarrow \zeta_0$ and $g_n^{(k)}(\zeta_n)=b; \ n=1,2,\ldots$  \\
Since $f_n(z)=b$ whenever $f_n^{(k)}(z)=b \Rightarrow f_n(z_n+\rho_n\zeta_n)=b$ and,\\

\qquad \qquad \qquad
$g_n(\zeta_n)=\rho_n^{-k}[f_n(z_n+\rho_n\zeta_n)-a]=\rho_n^{(k)}[b-a]
\rightarrow \infty$, this contradicts $$\lim_{n\rightarrow
\infty}g_n(\zeta_n)=g(\zeta_0)\neq \infty$$ So $g^{(k)}(\zeta) \neq
b$. Hence we get,

\begin{equation}\label{eq1}
 g^{(k)}(\zeta) =b+e^{A\zeta +B}
 \end{equation}

where $A$  and $B$ are two constants. We claim that $A=0$. Suppose that $A\neq0$; then

\begin{equation}\label{eq2}
g(\zeta)=\frac{b\zeta^{k}}{k!}+\frac{e^{A\zeta +B}}{A^{k}}+\frac{c_1\zeta^{k-1}}{(k-1)!}+\ldots +c_{k-1}\zeta +c_k
\end{equation}

where $c_1, c_2, \ldots,c_k$ are constants.
Let $g^{(k)}=a$. Then by\eqref{eq1},\eqref{eq2} and $g(\zeta)=0 \ \text{whenever} g^{(k)}(\zeta)=a \ $,
we have\\

\qquad \qquad \qquad \qquad  $\frac{b\zeta^k}{k!}+\frac{c_1\zeta^{k-1}}{(k-1)!}+\ldots +c_k+\frac{b-a}{A^k}=0$\\
This is a polynomial of degree k in $\zeta$ this polynomial has k solutions.which contradicts the fact that $g^{(k)}$ has infinitely many solutions. Thus we have,\\

\qquad \qquad \qquad \qquad \qquad \qquad \qquad $g^{(k)}(\zeta)=b+e^B$  \\

And\\

\qquad \qquad \qquad \qquad  $ g(\zeta)=(b+e^B)\frac{\zeta^k}{k!}+\frac{c_1\zeta^{k-1}}{(k-1)!}+\ldots +c_k$\\

Since $g$ is non constant, this contradicts $g(\zeta)=0\Leftrightarrow g^{(k)}(\zeta)=a$. Thus $\fr $ is normal in $D.$ This completes the proof of theorem.

\end{proof}

\section{Proof of Theorem 2.2}

\begin{proof}
Suppose that $\fr $ is not normal at $z_0\in\Delta$, then  by Lemma \ref{lem1}, there exist \\(a) a sequence of complex numbers $z_n \rightarrow z_0$, $|z_n|<r<1$\\
(b) a sequence of functions $f_n\in \mathcal F$ and\\
(c) a sequence of positive numbers $\rho_n \rightarrow 0$\\
such that $g_n(\zeta)=\rho_n^{-k}[f_n(z_n+\rho_n\zeta)-a(z_n+\rho_n\zeta)]$\\ converges locally uniformly to a non constant entire function $g$. Moreover $g$ is of order at most one.\\ Now we claim that \\
\qquad (a): \ $g(\zeta)=0\Leftrightarrow g^{(k)}(\zeta)=\varphi(z_0), $ \ where $\varphi(z)=\frac{a(z)-\al_1(z)a(z)-\al_0(z)a^{(k)}(z)}{\al_0(z)}$\\
\qquad (b): \ $g^{(k)}(\zeta)\neq B$, \ where $B=\frac{b(z_0)-\al_1(z_0)a(z_0)-\al_0(z_0)a^{(k)}(z_0)}{\al_0(z_0)}$ Note that $B$ is a constant.\\

Since

\begin{equation}\label{eq2.1}
g_n(\zeta)=\rho_n^{-k}[f_n(z_n+\rho_n\zeta)-a(z_n+\rho_n\zeta)]\rightarrow g(\zeta)
\end{equation}

we have

\begin{equation}\label{eq2.2}
g_n^{(k)}(\zeta)=f_n^{(k)}(z_n+\rho_n\zeta)-a^{(k)}(z_n+\rho_n\zeta)\rightarrow g^{(k)}(\zeta)
\end{equation}

Now suppose that $g(\zeta_0)=0.$ Then by Hurwitz's theorem, there exists $\zeta_n, \ \zeta_n\rightarrow \zeta_0$ such that\\

\qquad\qquad\qquad \qquad  $g_n(\zeta_n)=\rho_n^{-k}[f_n(z_n+\rho_n\zeta)-a(z_n+\rho_n\zeta)]=0$.\\

Thus\\

 \qquad\qquad\qquad \qquad $f_n(z_n+\rho_n\zeta)=a(z_n+\rho_n\zeta)$.\\

 Since$f(z)=a(z)\Leftrightarrow \al_0(z) f^{(k)}(z)+\al_1(z) f(z)=a(z)$, we have $\al_0(z) f^{(k)}(z)+\al_1(z) f(z)=a(z)$. \\

 Also\\

  $\frac{\al_0(z_n+\rho_n\zeta) f^{(k)}(z_n+\rho_n\zeta)+\al_1(z_n+\rho_n\zeta) f(z_n+\rho_n\zeta)}{\al_0(z_n+\rho_n\zeta)}=f_n^{(k)}(z_n+\rho_n\zeta)+ \frac{\al_1(z_n+\rho_n\zeta)}{\al_0(z_n+\rho_n\zeta)}f_n(z_n+\rho_n\zeta)$\\

  \qquad\qquad\qquad\qquad\qquad\qquad\qquad $=f_n^{(k)}(z_n+\rho_n\zeta)+ \frac{\al_1(z_n+\rho_n\zeta)}{\al_0(z_n+\rho_n\zeta)}[\rho_n g_n(\zeta)+a(z_n+\rho_n\zeta)]$\\
\begin{equation}\label{eq2.3}
\rightarrow g^{(k)}(\zeta)+a^{(k)}(z_0)+\frac{\al_1(z_0)}{\al_0(z_0)}a(z_0)
\end{equation}

  Therefore it follows that,\\

  \qquad\qquad $$g^{(k)}(\zeta_0)=\lim_{n\rightarrow\ity}[\frac{\al_0(z_n+\rho_n\zeta) f^{(k)}(z_n+\rho_n\zeta)+\al_1(z_n+\rho_n\zeta) f(z_n+\rho_n\zeta)}{\al_0(z_n+\rho_n\zeta)}]-a^{(k)}(z_0)-\frac{\al_1(z_0)}{\al_0(z_0)}a(z_0)$$\\
 $$=\lim_{n\rightarrow \ity}\frac{a(z_n+\rho_n\zeta)}{a_0(z_n+\rho_n\zeta)}-a^{(k)}(z_0)-\frac{\al_1(z_0)}{\al_0(z_0)}a(z_0)$$\\

 \qquad\qquad\qquad\qquad\qquad$=\frac{a(z_0)-\al_1(z_0)a(z_0)-\al_0(z_0)a^{(k)}(z_0)}{\al_0(z_0)}=\varphi(z_0).$\\

 Hence we have proved $g^{(k)}(\zeta)=\varphi(z_0)$ whenever $g(\zeta)=0$ \\

\qquad On the other hand, if $g^{(k)}(\zeta_0)=\varphi(z_0)$ then there exists $\zeta_n; \ \zeta_n\rightarrow \zeta_0, $  such that \\

\qquad\qquad\qquad  $f_n^{(k)}(z_n+\rho_n\zeta)-a^{(k)}(z_n+\rho_n\zeta)=\varphi(z_0)$\\

We have to show $$g(\zeta_0)=\lim_{n\rightarrow \ity}g_n(\zeta_n)=\lim_{n\rightarrow \ity}[f_n(z_n+\rho_n\zeta_n)-a(z_n+\rho_n\zeta_n)]=f(z_0)-a(z_0)=0$$\\

Now, assume  that $g^{(k)}(\zeta_0)=\varphi(z_0)$ by using assumption $(3)$ of the Theorem we get $f(z_0) -a(z_0) =0$, so is  $g(\zeta_0)=0.$ This shows that $g(\zeta)=0\Leftrightarrow g^{(k)}(\zeta)=\varphi(z_0) $\\

\qquad From \eqref{eq2.3}  we deduce that

\begin{equation}\notag
\frac{\al_0(z_n+\rho_n\zeta) f^{(k)}(z_n+\rho_n\zeta)+\al_1(z_n+\rho_n\zeta) f(z_n+\rho_n\zeta)-b(z_n+\rho_n\zeta)}{\al_0(z_n+\rho_n\zeta)}\rightarrow g^{(k)}(\zeta)+a^{(k)}(z_0)+\frac{\al_1(z_0)}{\al_0(z_0)}a(z_0)-\frac{b(z_0)}{\al_0(z_0)}
\end{equation}

\begin{equation}\label{eq2.4}
=g^{(k)}(\zeta)-\frac{b(z_0)-\al_0(z_0)a^{(k)}(z_0)-\al_1(z_0)a(z_0)}{\al_0(z_0)}= g^{(k)}(\zeta)-B
\end{equation}

Next we prove that $g^{(k)}(\zeta)\neq B.$ Suppose that there exists $\zeta_0$ satisfying $g^{(k)}(\zeta_0)=B.$ Then, by Hurwitz's theorem, there exists a sequence $\zeta_n; \ \zeta_n\rightarrow \zeta_0$  and by \eqref{eq2.4} \\

 $\{\al_0(z_n+\rho_n\zeta_n) f^{(k)}(z_n+\rho_n\zeta_n)+\al_1(z_n+\rho_n\zeta_n) f(z_n+\rho_n\zeta_n)\}-b(z_n+\rho_n\zeta_n)=0$\\

  From the assumption, we have $f_n(z_n+\rho_n\zeta_n)=b(z_n+\rho_n\zeta_n).$  Then we get\\

   $$g(\zeta_0)=\lim_{n\rightarrow\ity}\rho_n^{(k)}[f_n(z_n+\rho_n\zeta_n)-a(z_n+\rho_n\zeta_n)]$$\\

  \qquad $$=\lim_{n\rightarrow\ity}\rho_n^{(k)}[b(z_n+\rho_n\zeta_n)-a(z_n+\rho_n\zeta_n)]=\ity$$\\

  which is a contradiction. So $g^{(k)}(\zeta)\neq B.$\\

  \qquad Hence we get
\begin{equation}\label{eq2.5}
 g^{(k)}(\zeta)= B+e^{A\zeta+D}
 \end{equation}

  where $A$  and $D$ are two constants. We  claim that $A=0$. Suppose that $A\neq 0; $ then

\begin{equation}\label{eq2.6}
g(\zeta)=\frac{B\zeta^{k}}{k!}+\frac{e^{A\zeta +D}}{A^{k}}+\frac{c_1\zeta^{k-1}}{(k-1)!}+\ldots +c_{k-1}\zeta +c_k
\end{equation}
  where $c_1, c_2, \ldots,c_k$ are constants. Let $g^{(k)}(\zeta)=\varphi(z_0)$ then by \eqref{eq2.5}, \eqref{eq2.6} and $g(\zeta)=0\ \text{if and only if}\  g^{(k)}(\zeta)=\varphi(z_0)$\\

   So we get \\

   \qquad \qquad \qquad \qquad  $\frac{B\zeta^k}{k!}+\frac{c_1\zeta^{k-1}}{(k-1)!}+\ldots +c_k+\frac{B-\varphi(z_0)}{A^k}=0$\\

This is a polynomial of degree k in $\zeta$ this polynomial has k solutions.which contradicts the fact that $g^{(k)}$ has infinitely many solutions. Thus we have,\\

\qquad \qquad \qquad \qquad \qquad \qquad \qquad $g^{(k)}(\zeta)=B+e^D$  \\

And\\

\qquad \qquad \qquad \qquad $ g(\zeta)=(B+e^D)\frac{\zeta^k}{k!}+\frac{c_1\zeta^{k-1}}{(k-1)!}+\ldots +c_k$\\

Since $g$ is non constant, this contradicts $g(\zeta)=0\ \text{if and only if}\ g^{(k)}(\zeta)=\varphi(z_0)$. Thus $\fr $ is normal in $D.$ This completes the proof of theorem.

  \end{proof}

\section{Proof of Theorem 2.3}

\begin{proof}
Suppose $\fr$ is not normal in $\Delta$; without loss of generality we assume that $\fr$ is not normal at the point $z=0.$  Then by Lemma \ref{lem1}, there exist \\(a) a sequence of complex numbers $z_n \rightarrow 0$, $|z_n|<r<1$\\
(b) a sequence of functions $f_n\in \mathcal F$ and\\
(c) a sequence of positive numbers $\rho_n \rightarrow 0$\\
such that $\q$\\ converges locally uniformly to a non constant entire function $g$. Moreover $g$ is of order at most one.\\
Now we claim that   $g=0\text{ iff}\  g^{(k)}=a$ and $g^{(k+1)}=0 \  \text{whenever}\  g=0$\\ Let $g(\zeta_0)=0$. Then by Hurwitz's theorem there exist $\zeta_n;  \zeta_n\rightarrow \zeta_0$ such that\\

\qquad \qquad \qquad \qquad $g_n(\zeta_n)=\rho_n^{-k}[f_n(z_n+\rho_n\zeta_n)-a]=0 $\\

 Thus $f_n(z_n+\rho_n\zeta_n)=a$ since $f_n$ and $ f_n^{(k)}$ share $z=a$ IM, we have\\

\qquad \qquad \qquad \qquad \qquad \qquad  $g_n^{k}(\zeta_n)=f_n^{(k)}(z_n+\rho_n\zeta_n)=a$\\

and\\

\qquad \qquad \qquad \qquad \qquad \qquad  $g_n^{(k+1)}(\zeta_n)=\rho_nf_n^{(k+1)}(z_n+\rho_n\zeta_n)$\\

  which implies that\\

  $$g^{(k)}(\zeta_0)=\lim_{n \rightarrow \infty}g_n^{(k)}(\zeta_n)=a$$\\ and $$g^{(k+1)}(\zeta_0)=\lim_{n \rightarrow \infty}g_n^{(k+1)}(\zeta_n)=0$$\\

  Thus we get, $g^{(k)}=a \text{ whenever }  g=0$  and $g^{(k+1)}=0  \text{ whenever }  g=0$.\\

\qquad On then other hand, if $g^{(k)}(\zeta_0)=a$ then there exit $\zeta_n \rightarrow \zeta_0$ such that \\
$g_n^{(k)}(\zeta_n)=f_n^{(k)}(z_n+\rho_n\zeta_n)=a,\  n=1,2,\ldots$ hence $f_n(z_n+\rho_n\zeta_n)=a$ and $g_n(\zeta_n)=0$ for n=1,2,\ldots  thus $$g(\zeta_0)=\lim_{n\rightarrow \infty}g_n(\zeta_n)=0.$$ This shows that $g=0$ whenever
$g^{(k)}=a$.\\

 Hence $g=0$ if and only if  $g^{k}=a$ and $g^{(k+1)}=0 \text{ whenever } g=0$.\\

Now using Lemma \ref{lem2} and Nevanlinna's first fundamental theorem, we have\\

\qquad \qquad  $T(r,g)\leq \overline N(r,g) +N(r,\frac{1}{g})+N(r,\frac{1}{g^{(k)}-a})-N(r,\frac{1}{g^{(k+1)}})+S(r,g)$\\

\qquad \qquad \qquad \qquad = $N(r,\frac{1}{g})+N(r,\frac{1}{g^{(k)}-a}) -N(r,\frac{1}{g^{(k+1)}})+S(r,g)$\\

\qquad \qquad \qquad \qquad $\leq N(r,\frac{1}{g^{(k)}-a})-\overline N(r,\frac{1}{g^{(k+1)}})+S(r,g)$\\

\qquad \qquad \qquad \qquad $\leq T(r,\frac{1}{g^{(k)}-a})-\overline N(r,\frac{1}{g^{(k+1)}})+S(r,g)$\\

\qquad \qquad \qquad \qquad$\leq T(r,g^{(k)}-a)-\overline N(r,\frac{1}{g^{(k+1)}})+S(r,g)$\\

\begin{equation}\label{eq3.1}
\leq T(r,g)-\overline N(r,\frac{1}{g^{(k+1)}})+S(r,g)
\end{equation}

 Thus we get\\

\begin{equation}\label{eq3.2}
\overline N(r,\frac{1}{g^{(k+1)}})=S(r,g)
\end{equation}

 by \eqref{eq3.1}, \eqref{eq3.2} and the claim($g=0$ if and only if $g^{(k)}=a, g^{(k)}=g^{(k+1)}=0 \text{ whenever} \ g=0$) we get a contradiction: $T(r,g)=S(r,g).$\\ Hence the theorem.

\end{proof}


\begin{thebibliography}{00}

\bibitem{Hay} W.K. Hayman, Meromorphic Functions, Claredon Press, Oxford, 1964.
\bibitem{Schiff} J. Schiff,  Normal Families, Springer-Verlag, Berlin, 1993.
\bibitem{Yang}L. Yang, Value Distribution Theory, Springer- Verlag,Berlin, 1993.
\bibitem{Zalc}L. Zalcman, Normal Families: New perspective, Bulletin of American Mathematical Society \textbf{35} (1998), 215-230.
\bibitem{Zalc 1}L. Zalcman, A heuristic principle in complex function theory, The American Mathematical Momthly \textbf{82} (1975),813-817.
\bibitem{CH}H.H. Chen and X.H. Hua, Normal Values Concerning Shared Values, Israel Journal of Mathematics  \textbf{115} (2000), 355-362.
\bibitem{FX}Ming-Ling Fang and Yan Xu, Normal Families of Holomorphic Functions and Shared Values, Israel Journal of Mathematics   \textbf{129} (2002), 125-141.
\bibitem{Sch}W. Schwick, Sharing Values and Normality, Archiv der Mathematik\textbf{59},(1992),50-54.
\bibitem{MS}E. Mues and N. Steinmentz, Meromorphe Functionen, die mit ihrer Ableitung Werte teilen, Manuscripta Mathematica \textbf{29} (1979), 195-206.
\bibitem{MLF}Ming-Ling Fang, Picard Values and Normality criterion, Bull. Korean Math. Soc. \textbf{38} (2001), No. 2, 379-387.
\bibitem{YX}Yan Xu, Normality Criterion concerning Sharing Functions, Houstan Journal of Mathematics \textbf{32}, No. 3, (2006) 945-954.
\bibitem{XX}Jiying Xia, Yan Xu, Normality Criterion concerning Sharing Functions II, Bull. Malays. Math. Sci. Soc.(2) \textbf{33} (3)(2010), 479-486.
\bibitem{CF}J.F. Chen, M.L. Fang, Normal Families And Shared Functions of Meromorphic Functions, Israel Journal of Mathematics  \textbf{180} (2010), 129-142.
\bibitem{JFC}Jun-Fan Chen, Normal Families and Shared Sets of Meromorphic Functions, Rocky Mountain Journal of Mathematics, \textbf{41}, No. 1, (2011).






\end{thebibliography}
\end{document}